\documentclass[11pt]{article}
\usepackage{amsmath,amsthm,amsfonts}
\newcommand\pd[2]{\frac{\partial#1}{\partial#2}}
\renewcommand{\=}{\doteq}

\newcommand{\X}{\mathfrak{X}}

\newtheorem{thm}{Theorem}[section]
 \newtheorem{prop}[thm]{Proposition}
 \newtheorem{lemma}[thm]{Lemma}
\theoremstyle{definition}

\theoremstyle{definition}
 \newtheorem{rem}[thm]{Remark}
\numberwithin{equation}{section}
\numberwithin{equation}{section}
\begin{document}
\title{\bf  Moufang symmetry XII.\\
Reductivity and hidden associativity of\\
infinitesimal Moufang transformations
}
\author{Eugen Paal}
\date{}
\maketitle
\thispagestyle{empty}
\begin{abstract}
It is shown how integrability of the generalized Lie equations of continous Moufang transformatiosn is related to the reductivity conditions and Sagle-Yamaguti identity. 
\end{abstract}

\section{Introduction}

In this paper we proceed explaing the Moufang symmetry. It is shown how integrability of the generalized Lie equations of a local analytic Moufang loop is related to the reductivity conditions and Sagle-Yamaguti identity. The paper can be seen as a continuation of \cite{Paal1,Paal2,Paal3,Paal11,Paal_AAM}.

\section{Generalized Lie equations}

In \cite{Paal1} the \emph{generalized Lie equations} (GLE) of a local analytic Moufang loop $G$ were found. These read
\begin{subequations}
\label{gle_S}
\begin{align}
u^{s}_{j}(g)\pd{(S_gA)^{\mu}}{g^{s}}+T^{\nu}_{j}(A)\pd{(S_gA)^{i}}{A^{\nu}}+P^{\nu}_{j}(S_gA)&=0\\
v^{s}_{j}(g)\pd{(S_gA)^{\mu}}{g^{s}}+P^{\nu}_{j}(h)\pd{(S_gA)^{i}}{A^{\nu}}+T^{\nu}_{j}(S_gA)&=0\\
w^{s}_{j}(g)\pd{(S_gA)^{\mu}}{g^{s}}+S^{\nu}_{j}(h)\pd{(S_gA)^{i}}{A^{\nu}}+S^{\nu}_{j}(S_gA)&=0
\end{align}
\end{subequations}
where $gh$ is the product of $g$ and $h$, and the auxiliary functions $u^s_j$, $v^s_j$, $w^s_j$ and
$S^\mu_j$, $T^\mu_j$, $P^\mu_j(g)$ are related with the  constraints
\begin{gather}
u^s_j(g)+v^s_j(g)+w^s_j(g)=0\\
S^\mu_j(A)+T^\mu_j(A)+P^\mu_j(A)=0
\end{gather}
For $T_gA$ the GLE read
\begin{subequations}
\label{gle_T}
\begin{align}
v^{s}_{j}(g)\pd{(T_gA)^{\mu}}{g^{s}}+S^{\nu}_{j}(A)\pd{(T_gA)^{i}}{A^{\nu}}+P^{\nu}_{j}(T_gA)&=0\\
u^{s}_{j}(g)\pd{(T_gA)^{\mu}}{g^{s}}+P^{\nu}_{j}(h)\pd{(T_gA)^{i}}{A^{\nu}}+S^{\nu}_{j}(T_gA)&=0\\
w^{s}_{j}(g)\pd{(T_gA)^{\mu}}{g^{s}}+T^{\nu}_{j}(h)\pd{(T_gA)^{i}}{A^{\nu}}+T^{\nu}_{j}(T_gA)&=0
\end{align}
\end{subequations}
In this paper we inquire integrability of GLE (\ref{gle_S}a--c). and (\ref{gle_T}a--c) 
Triality \cite{Paal2} considerations are wery helpful.

\section{Generalized Maurer-Cartan equations and Yamagutian}

Recall from \cite{Paal_AAM} that for $x$ in  $T_e(G)$ the infinitesimal translations of $G$ are defined by
\begin{equation*}
S_x\=x^j S^\nu_j(A)\pd{}{A^\nu},\quad
T_x\=x^j T^\nu_j(A)\pd{}{A^\nu},\quad
P_x\=x^j P^\nu_j(A)\pd{}{A^\nu}\quad \in T_A(\X)
\end{equation*}
with constriant
\begin{equation*}
S_x+T_x+P_x=0
\end{equation*}
Following triality \cite{Paal2} define the Yamagutian $Y(x;y)$ by
\begin{equation*}
6Y(x;y)=[S_x,S_y]+[T_x,T_y]+[P_x,P_y]
\end{equation*}
We know  from \cite{Paal_AAM} the generalized Maurer-Cartan equations:
\begin{subequations}
\label{m-c_ST}
\begin{align} 
[S_{x},S_{y}]&=S_{[x,y]}-2[S_{x},T_{y}]\\
[T_{x},T_{y}]&=T_{[y,x]}-2[T_{x},S_{y}]\\
[S_{x},T_{y}]&=[T_{x},S_{y}],\quad \forall x,y \in T_e(G) 
\end{align}
\end{subequations}
The latter can be written \cite{Paal2} as follows:
\begin{subequations}
\label{lr-yam_ST}
\begin{align}
[S_{x},S_{y}]&=2Y(x;y)+\frac{1}{3}S_{[x,y]}+\frac{2}{3}T_{[x,y]}\\
[S_{x},T_{y}]&=-Y(x;y)+\frac{1}{3}S_{[x,y]}-\frac{1}{3}T_{[x,y]}\\
[T_{x},T_{y}]&=2Y(x;y)-\frac{2}{3}S_{[x,y]}-\frac{1}{3}T_{[x,y]}
\end{align}
\end{subequations}

\section{Reductivity}

Define the  (secondary) auxiliary functions of $G$ by
\begin{align*}
S^\mu_{jk}(A)
&\=S^\nu_k(A)\pd{S^\mu_j(A)}{A^\nu}-S^\nu_j(g)\pd{S^\mu_k(A)}{A^\nu}\\
T^\mu_{jk}(A)
&\=T^\nu_k(A)\pd{T^\mu_j(A)}{A^\nu}-T^\nu_j(g)\pd{T^\mu_k(A)}{A^\nu}\\
P^\mu_{jk}(A)
&\=P^\nu_k(A)\pd{P^\mu_j(A)}{A^\nu}-P^\nu_j(g)\pd{P^\mu_k(A)}{A^\nu}
\end{align*}
The Yamaguti functions $Y^\mu_{jk}$ are defined by
\begin{equation*}
6Y^\mu_{jk}(A)\=S^\mu_{jk}(A)+T^\mu_{jk}(A)+P^s_{jk}(A)
\end{equation*}
In \cite{Paal3} we proved
\begin{thm}
The integrability conditons of the GLE (\ref{gle_S}a--c) (\ref{gle_T}a--c) read, respectively,
\begin{subequations}
\label{gle2yam_ST}
\begin{align}
Y^s_{jk}(g)\pd{(S_gA)^\mu}{g^s}+Y^\nu_{jk}(A)\pd{(S_gA)^\mu}{A^\nu}&=Y^\mu_{jk}(S_gA)\\
Y^s_{jk}(g)\pd{(T_gA)^\mu}{g^s}+Y^\nu_{jk}(A)\pd{(T_gA)^\mu}{A^\nu}&=Y^\mu_{jk}(T_gA)
\end{align}
\end{subequations}
\end{thm}
Consider the first-order approximation of the integrability conditions (\ref{gle2yam_ST}a) and (\ref{gle2yam_ST}b). We need 
\begin{lemma}
One has
\begin{equation}
\label{assoc2yam_bir}
Y^\mu_{jk}=l^\mu_{jk}+\frac{1}{3}C^s_{jk}(S^\mu_s-T^\mu_s)
\end{equation}
\end{lemma}

\begin{proof}
Use formula (\ref{lr-yam_ST}b).
\end{proof}

Introduce the Yamaguti constants $Y^i_{jkl}$ by
\begin{equation*}
Y^i_{jk}(g)=Y^i_{jkl}g^l+O(g^2)
\end{equation*} 
Then, by defining \cite{Paal1} the third-order associators $l^i_{jkl}$ by 
\begin{equation*}
l^i_{jk}(g)=l^i_{jkl}g^l+O(g^2)
\end{equation*}
it follows from Lemma \ref{assoc2yam_bir} that
\begin{equation}
\label{yam3}
Y^i_{jkl}=l^i_{jkl}+\frac{1}{3}C^s_{jk}C^i_{sl}
\end{equation}
Now we can calculate:
\begin{gather*}
Y^i_{jk}(S_gA)
=Y^\mu_{jk}(A)+\pd{Y^\mu_{jk}(A)}{A^\nu}S^\nu_l(h)g^l+O(g^2)\\
Y^s_{jk}(g)\pd{(S_gA)^\mu}{g^s}+Y^\nu_{jk}(A)\pd{(S_gA)^\mu}{A^\nu}
=Y^s_{jkl}g^l S^\mu_s(h)+Y^\nu_{jk}(A)\left(\delta^\mu_\nu+\pd{S^\mu_{l}(A)}{A^\nu}g^l\right)+O(g^2)
\end{gather*}
Substituting the latter into (\ref{gle2yam_ST}a) and compare the coefficients at $g^l$ and replace.
By repeating these calculations for (\ref{gle2yam_ST}b)  we obtain the \emph{reductiovity conditions}
\begin{subequations}
\label{pre-red_bir}
\begin{align}
S^\nu_l(g)\pd{Y^\mu_{jk}(A)}{A^\nu}-Y^\nu_{jk}(A)\pd{S^\mu_{l}(A)}{A^\nu}&=Y^s_{jkl}S^\mu_s(A)\\
T^\nu_l(g)\pd{Y^\mu_{jk}(A)}{A^\nu}-Y^\nu_{jk}(A)\pd{T^\mu_{l}(A)}{A^\nu}&=Y^s_{jkl}T^\mu_s(A)
\end{align}
\end{subequations}
Let us rewrite these differential equations as commutation relations.

In the tangent algebra $\Gamma$ of $G$ define the the ternary \emph{Yamaguti brackets} \cite{Yam63} $[\cdot,\cdot,\cdot]$ by
\begin{equation*}
[x,y,z]^i\=6Y^i_{jkl}x^jy^kz^l
\end{equation*}
Multiply (\ref{yam3}) by $6x^jy^kz^l$. Then we have
\begin{align*}
[x,y,z]
&=6(x,y,z)+2[[x,y],z]\\
&=[x[y,z]]-[y[x,z]]+[[x,y],z]
\end{align*}
Now from (\ref{pre-red_bir}) it is easy to infer 
\begin{thm}[reductivity]
The infinitesimal Moufang transformations satisfy the reductivity conditions
\begin{subequations}
\label{red_bir}
\begin{align}
6[Y(x;y),S_z]&=S_{[x,y,z]}\\
6[Y(x;y),T_z]&=T_{[x,y,z]}\\
6[Y(x;y),P_z]&=P_{[x,y,z]}
\end{align}
\end{subequations}
\end{thm}

\begin{proof}
Commutation relations (\ref{red_bir}a,b) are evident from  (\ref{pre-red_bir}a,b) and  (\ref{red_bir}c) easily follows by adding (\ref{red_bir}a) and  (\ref{red_bir}b).
\end{proof}

\section{Sagle-Yamaguti identity and hidden associativity}

Define the triality conjugated translations
\begin{equation*}
P^+\=S-T,\quad 
S^+\=T-P,\quad
T^+\=P-S
\end{equation*}
One can easily see the inverse conjugation:
\begin{equation*}
3P\=T^+-S^+,\quad 
3T\=S^+-T^+,\quad
3P\=T^+-S^+
\end{equation*}

\begin{thm}[reductivity]
The infinitesimal Moufang transformations satisfy the reductivity conditions
\begin{subequations}
\label{red+}
\begin{align}
6[Y(x;y),S^+_z]&=S^+_{[x,y,z]}\\
6[Y(x;y),T^+_z]&=T^+_{[x,y,z]}\\
6[Y(x;y),P^+_z]&=P^+_{[x,y,z]}
\end{align}
\end{subequations}
\end{thm}

\begin{proof}
Evident corollary from formulae (\ref{red_bir}).
\end{proof}

From \cite{Paal2} we know

\begin{prop}
Let $(S,T)$ be a Moufang-Mal'tsev pair. Then
\begin{subequations}
\label{lrm+}
\begin{align}
6Y(x;y)
&=[P^{+}_{x},P^{+}_{y}]+P^{+}_{[x,y]}\\
&=[T^{+}_{x},T^{+}_{y}]+T^{+}_{[x,y]}\\
&=[S^{+}_{x},S^{+}_{y}]+S^{+}_{[x,y]}
\end{align}
\end{subequations}
for all $x,y$ in $M$.
\end{prop}

\begin{thm}[hidden associativity]
The Yamagutian $Y$ of $(S,T)$ obey the commutation relations
\begin{equation}
\label{yam_com_bir}
6[Y(x;y),Y(z,w)]=Y([x,y,x],w)+Y(z;[x,y,w])
\end{equation}
if the following Sagle-Yamaguti identity \cite{Yam62,Yam63} holds:
\begin{equation}
\label{sagle-yamaguti}
[x,y,[z,w]]=[[x,y,z],w]+[z,[x,y,w]]
\end{equation}
\end{thm}

\begin{proof}
We calculate the Lie bracket $[Y(x;y),Y(z,w)]$ from the Jacobi identity
\begin{equation}
\label{jacobi-temp_bir}
[[Y(x;y),S^+_z],S^+_w]+[[S^+_z,S^+_w],Y(x;y]+[[S^+_w,Y(x;y),S^+_z]=0
\end{equation}
and formulae (\ref{lrm+}). We have
\begin{align*}
6[[Y(x;y),S^+_z],S^+_w]
&=[S^+_{[x,y,z]},S_w]\\
&=6Y([x,y,z];w)-S^+_{[[x,y,z],w]}\\
6[[S^+_z,S^+_w],Y(x;y]
&=36[Y(z;w),Y(x,y)]-6[S^+_{[z,w]},Y(x;y)]\\
&=36[Y(z;w),Y(x,y)]-S^+_{[x,y,[z,w]]}\\
6[[S^+_w,Y(x;y),S^+_z]
&=6Y(z;[x,y,w])-S^+_{[z,[x,y,w]]}
\end{align*}
By substituting these relations into (\ref{jacobi-temp_bir}) we obtain
\begin{equation*}
36[Y(x;y),Y(z,w)]-6Y([x,y,x],w)-6Y(z;[x,y,w])
=S^+_{[x,y,[z,w]]-[[x,y,z],w]-[z,[x,y,w]]}
\end{equation*}
The latter realtion has to be triality invariant. This mean that
\begin{subequations}
\label{stp+bir}
\begin{align}
S^+_a
&=T^+_a=P^+_a\\
&=36[Y(x;y),Y(z,w)]-6Y([x,y,x],w)-6Y(z;[x,y,w])
\end{align}
\end{subequations}
where
\begin{equation*}
a=[x,y,[z,w]]-[[x,y,z],w]-[z,[x,y,w]]
\end{equation*}
But it easily follows from (\ref{stp+bir}a) that
\begin{equation*}
S_a=T_a=P_a=0
\end{equation*}
and due to $a=0$ commutation relations (\ref{yam_com_bir}) hold.
\end{proof}

\begin{rem}
A.~Sagle \cite{Sagle} and K.~Yamaguti proved \cite{Yam62} that the identity (\ref{sagle-yamaguti}) is equivalent to the Mal'tsev identity. In terms of Yamaguti \cite{Yam63} one can say that the Yamagutian $Y$ is a \emph{generalized representation} of the (tangent) Mal'tsev algebra $\Gamma$ of $G$.
\end{rem}

\section*{Acknowledgement} 

Research was in part supported by the Estonian Science Foundation, Grant 6912.

\bigskip\noindent
Department of Mathematics\\
Tallinn University of Technology\\
Ehitajate tee 5, 19086 Tallinn, Estonia\\ 
E-mail: eugen.paal@ttu.ee

\end{document}